\documentclass[draft,10pt]{article}

\usepackage{psfig}
\usepackage{amssymb,amsmath,amsthm,amsfonts}

\newtheoremstyle{theorem}
  {}
  {}
  {\itshape}
  {}
  {\bfseries}
  {.}
  {.5em}
  {}

\newtheoremstyle{lemma}
  {}
  {}
  {\itshape}
  {}
  {\ttfamily}
  {.}
  {.5em}
  {}

\newtheoremstyle{plaintext}
  {}
  {}
  {\upshape}
  {}
  {\ttfamily}
  {.}
  {.5em}
  {}

\theoremstyle{theorem}
\newtheorem{theorem}{Theorem}

\theoremstyle{lemma}
\newtheorem{lemma}{Lemma}

\newtheorem{corollary}{Corollary}

\theoremstyle{plaintext}

\newtheorem{definition}{Definition}
\newtheorem{remark}{Remark}

\newcounter{savetheorem1}
\newcounter{savetheorem2}

\def\fg{\leqslant_{\text{f.g.}}}
\def\rank{\mbox{rank}}

\begin{document}

\title{On the Hanna Neumann Conjecture}

\author{Toshiaki Jitsukawa \and Bilal Khan \and Alexei G. Myasnikov}
\date{}

\maketitle

\begin{abstract}
The {\bf Hanna Neumann conjecture} states that if $F$ is a free group,
then for all nontrivial finitely generated subgroups $H,K \leqslant
F$,
\begin{eqnarray}
\rank(H \cap K)-1 \leqslant \left[\rank(H)-1\right] \left[\rank(K)-1\right] \nonumber
\end{eqnarray}
Where most papers to date have considered a direct graph theoretic
interpretation of the conjecture, here we consider the use of
monomorphisms.  We illustrate the effectiveness of this approach with
two results.  First, we show that for any finitely generated groups $H,
K \leqslant F$ either the pair $H,K$ or the pair $H^{-},K$ satisfy the
Hanna Neumann conjecture; here ${}^{-}$ denotes the automorphism which
sends each generator of $F$ to its inverse.  Next, using particular
monomorphisms from $F$ to $F_2$, we obtain that if the Hanna Neumann
conjecture is false then there is a counterexample $H, K \leq F_2$
having the additional property that all the branch vertices in the
foldings of $H$ and $K$ are of degree $3$, and all degree $3$ vertices
have the same local structure or ``type''.
\end{abstract}

\section{Introduction}

H. Neumann proved in \cite{HNeumann} that any nontrivial subgroups $H,K
\fg F$ (finitely generated) must satisfy
\begin{eqnarray}
\label{hntheorem}
 \; rank(H\cap K) - 1 \leqslant 2 [rank(H) - 1]
[rank(K) - 1],
\end{eqnarray}
and so improved Howson's earlier result \cite{Howson} that $H\cap K$ is
finitely generated. The stronger assertion obtained by omitting the
factor of $2$ in (\ref{hntheorem}) has come to be known as the Hanna
Neumann conjecture.  In \cite{Burns}, R. Burns improved H. Neumann's
bound by showing that
\begin{eqnarray*}
rank(H\cap K) - 1 & \leqslant & 2 [rank(H) - 1] [rank(K) - 1] \\
                  &      & - \min(rank(H) - 1 , rank(K) - 1).
\end{eqnarray*}

In 1983, J. Stallings introduced the notion of a {\em folding} and
showed how to apply these objects to the study of subgroups of free
groups \cite{Stallings}.  Stallings's approach was applied by S.
Gersten in \cite{Gersten} to solve certain special cases of the
conjecture, and similar techniques were developed over a sequence of
papers by W. Imrich \cite{Imrich2,Imrich}, P. Nickolas \cite{Nickolas},
and B. Servatius \cite{Servatius} who gave alternate proofs of Burns'
bound and resolved special cases of the conjecture.  In 1989, W.
Neumann showed that the conjecture is true ``with probability 1'' for
randomly chosen subgroups of free groups \cite{WNeumann}, and proposed
a stronger form of the conjecture.  In 1992, G. Tardos proved in
\cite{Tardos92} that the conjecture is true if one of the two subgroups
has rank 2.  Then, in 1994, W. Dicks showed that the strong Hanna
Neumann conjecture is equivalent to a conjecture on bipartite graphs,
which he termed the Amalgamated Graph conjecture \cite{Dicks}.  In
1996, G. Tardos used Dicks' method to give the first new bound for the
general case in \cite{Tardos96}, where he proved that for any  $H, K
 \leq F$ with $rank(H), rank(K) \geq 3$,
\begin{eqnarray*}
rank(H\cap K) - 1 & \leqslant & 2 [rank(H) - 1] [rank(K) - 1] \\
                  &      & - [rank(H)-1] - [rank(K)-1]
\end{eqnarray*}
 Since then, W. Dicks and E.Formanek \cite{rank3} proved that
\begin{eqnarray*}
rank(H \cap K)-1 & \leq & [rank(H)-1] [rank(K) - 1] + \\
 & & \max\{rank(H)-3, 0\} \max \{rank(K)-3, 0\},
\end{eqnarray*}
This resolved the conjecture for the case when one of the subgroups
has rank at most $3$.

The conjecture was also recently solved in the special case when
one of the two groups, say $H$, has a generating set consisting of
positive words (i.e. a set of words in which no generator of $F$
has negative exponent).  Specifically, it was shown by J. Meakin
and P. Weil \cite{MeakinWeil}, and independently by B. Khan
\cite{KhanHNC} that if there is some automorphism of $F$ which
carries a generating set of $H$ to a set of positive words, then
the conjecture holds for $H$ and any nontrivial $K \fg F$.

\vspace{0.2in}

Recall that an automorphism $\sigma$ of $F(X)$ is called {\em
length-preserving} if $\forall u\in F$, $|u^\sigma| = |u|$, i.e.
$(X^\pm)^\sigma = X^\pm$ where $X^\pm = X \cup X^{-1}$.  In
section \ref{results}, we shall prove the following two theorems:

\begin{theorem}
\label{first-theorem}
\setcounter{savetheorem1}{\value{theorem}}
\addtocounter{savetheorem1}{-1}
Let $F_2=F(a,b)$.  Take $\sigma \in Aut(F_2)$ to be any
length-preserving automorphism having no non-trivial fixed points, and
let $\tau$ be any monomorphism
\begin{eqnarray*}
\tau : \left.
\begin{array}{rcl}
    a & \mapsto & a  w_a  a \\
    b & \mapsto & b  w_b  b,
\end{array}
\right.
\end{eqnarray*}
where $w_a, w_b \in F_2$ are arbitrary elements for which the words $a
w_a a$ and $b w_b b$ are reduced as written. Then for all nontrivial
$H,K \fg F_2$, either the pair $H, K$ or the pair $H^{\tau\sigma}, K$
satisfy the conjecture.
\end{theorem}

\begin{theorem}
\label{first-theorem-part2}
\setcounter{savetheorem2}{\value{theorem}}
\addtocounter{savetheorem2}{-1} Let $F=F(X)$ and ${}^{-}$ be the
automorphism given by $x \mapsto x^{-1}$ (for each $x\in X$).  Then
for all nontrivial $H,K \fg F$, either the pair $H, K$ or the pair
$H^{-}, K$ satisfy the conjecture.
\end{theorem}

\vspace{0.2in}

\begin{figure}[h]
\centering{\mbox{\psfig{figure=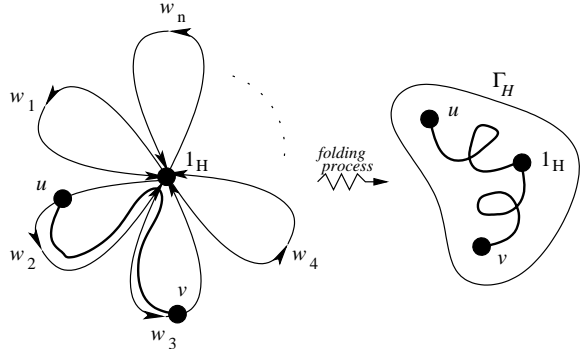}}}
\caption{Contructing a folding from a rose.}
\label{folding}
\end{figure}

Recall that given $H=\langle w_1, \cdots, w_n \rangle \leq F$, one may
determine the associated Stallings' folding $\Gamma_{H}=(V_H,E_H)$, by
the following constructive procedure (see \cite{Stallings}): Construct
$n$ directed cycles $c_1=(V_1,E_1)$, $\ldots$, $c_n=(V_n,E_n)$, where
$|V_i|=|w_i|$.  Then pick one vertex from each of the cycles, and
identify this subset of vertices, denoting the resulting vertex $1_H$.
Label the edges of cycle $c_i$ by successive letters of $w_i$, starting
at vertex $1_H$.  Finally, repeatedly identify pairs of edges $e,e'$
for which
\[
label(e)=label(e') \;\wedge\; [head(e)=head(e') \vee tail(e)=tail(e')] \nonumber.
\]
Each such identification is called an {\em edge-folding} and we say
that the edge $e$ (as well as $e'$) was {\em folded}.  Figure
\ref{folding} illustrates the process, which terminates in finitely
many steps yielding the folding $\Gamma_H$. It is easy to verify that
the folding so obtained is well-defined, and moreover, is independent
of the choice of generating set for $H$. It is not hard to see that
the rank of $H$ is precisely $|E_H| - |V_H| + 1$.

\begin{figure}[htbp]
\centering{\mbox{\psfig{figure=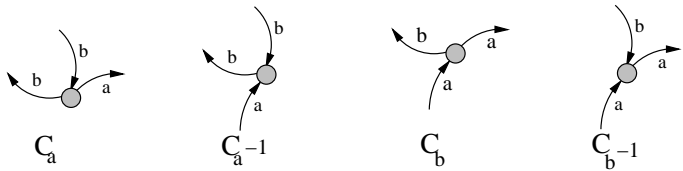}}}
\caption{Local structure of vertex labels in a folding.}
\label{local-types}
\end{figure}

Now we consider subgroups of $F_2$:  If $H\fg F_2$ then $\Gamma_H$
has vertices of undirected degree $\leqslant 4$, where by
``undirected degree'' $d=d_H(v)$ of a vertex $v$ we mean the sum
of the number of outgoing and incoming edges at $v$.  Put
$d_i(\Gamma_H) = \left| \left\{ v \in V_H | \text{ d}_H(v) =
i\right\} \right|$, for $i=1,2,3,4$. Vertices of degree $3$ may be
classified into $4$ types, denoted $C_a, C_b, C_{a^{-1}},
C_{b^{-1}}$, based on the labels of the incident edges (see figure
\ref{local-types}). For each $x \in \{a^\pm,b^\pm\}$, we define
$C_x(\Gamma_H)$ to be the number of degree $3$ vertices of type
$C_x$ in $\Gamma_H$.  The rank of $H$ can be computed by the
formula
\[
rank(H) = d_4(\Gamma_H) + \frac{d_3(\Gamma_H)}{2} -
\frac{d_1(\Gamma_H)}{2} + 1
\]

\vspace{0.2in}

The graph-theoretic approach to Hanna Neumann's conjecture is
based on the following key observation \cite{Stallings,Ilya}.
Consider the product automaton $\Gamma_H \times \Gamma_K$, whose
vertex set is $V_H\times V_K$ and two vertices $(u_1, v_1)$ and $
(u_2, v_2)$ are connected by an edge labelled $x$ iff both $(u_1,
v_1) \in E_H$ and $(u_2, v_2)\in E_K$ have label $x$. Consider the
connected component $\Delta$ of $\Gamma_H \times \Gamma_K$ which
contains $(1_H, 1_K)$.  The core of $\Delta$ is obtained by
repeated deletion of all vertices of degree $1$ different from
$(1_H, 1_K)$.  It is not hard to see that the core of $\Delta$ is
precisely $\Gamma_{H\cap K}$.

Many proofs of the Hanna Neumann conjecture (for groups of particular
ranks) require a case-by-case analysis based on the numbers and types
of degree~$3$ vertices present in the foldings of $H,K$.  The next
theorem has implications on the number of cases which need to be
considered in such arguments; it is proved in section \ref{results}.

\begin{theorem}
\label{second-theorem}
\setcounter{savetheorem2}{\value{theorem}}
\addtocounter{savetheorem2}{-1}
Let $F=F(X)$ be a non-abelian free group.  There is a monomorphism
$\phi_0:F\rightarrow F_2$ into the free group of rank $2$ such that:
for any groups $H, K \leqslant F$, the foldings of $H^{\phi_0},
K^{\phi_0} \fg F_2$ have the property that all their branch vertices
are of degree $3$, and all degree $3$ vertices have the same type.
\end{theorem}

The previous theorem has the following immediate corollary:

\begin{corollary}
\label{third-theorem}
\setcounter{savetheorem2}{\value{theorem}}
\addtocounter{savetheorem2}{-1}
If the Hanna Neumann conjecture is false then there is a
counterexample $H, K \fg F_2$ having the additional property that all
the branch vertices in the foldings of $H$ and $K$ are of degree $3$,
and all degree $3$ vertices have the same type.
\end{corollary}

\section{Preliminaries}

The numbers $C_x(\Gamma_H)$ and $C_x(\Gamma_K)$ allow one to
compute upper bounds on the numbers of vertices of degree $3$ in
$\Gamma_H \times \Gamma_K$ and hence in $\Gamma_{H\cap K}$.
Observe that by considering a suitable conjugate of $H$ and $K$ we
can always assume $d_1(\Gamma_H) = d_1(\Gamma_K) = 0$.
Furthermore,
\begin{eqnarray*}
d_4(\Gamma_{H\cap K}) & \leqslant & d_4(\Gamma_H \times \Gamma_K)
= d_4(\Gamma_H) d_4(\Gamma_K), \text{ and}\\
d_3(\Gamma_{H\cap K}) & \leqslant&  d_3(\Gamma_H \times \Gamma_K)
\\
& = & d_4(\Gamma_H) d_3(\Gamma_K) + d_3(\Gamma_H) d_4(\Gamma_K) +
\sum_{x \in \{a,b\}^\pm} C_x(\Gamma_H) C_x(\Gamma_K)
\end{eqnarray*}

\begin{definition}
Given two subgroups $H,K \fg F_2$ and $x~ \in~\{a, b\}^\pm$, we
define
\begin{eqnarray*}
\delta_x(H,K) & = & \min \left\{ \frac{C_x(\Gamma_H)}{d_3(\Gamma_H)}, \frac{C_x(\Gamma_K)}{d_3(\Gamma_K)} \right\} \\[0.1in]
\delta(H,K)   & = &\max_{x \in \{a^\pm,b^\pm\}} \delta_x(H,K)
\end{eqnarray*}
and put $\mu(H,K)$ to be any $x \in \{a^\pm,b^\pm\}$ for which
$\delta_x(H,K) = \delta(H,K)$.
\end{definition}

\begin{remark}
\label{walter-neumann} Walter Neumann \cite{WNeumann} showed that
if $H,K \fg F_2$ are a counterexample to the conjecture, then
$\delta(H,K) > \frac{1}{2}$. We outline his argument here, in
graph-theoretic notation. Using a simple and beautiful argument
from convexity theory, he showed that if $\delta(H,K) \leqslant
\frac{1}{2}$ then
\[
 \sum_{x \in \{a,b\}^\pm} C_x(\Gamma_H) C_x(\Gamma_K) \leqslant
 \frac{1}{2} d_3(\Gamma_H) d_3(\Gamma_K).
 \]
 It follows then that
\begin{eqnarray*}
\rank(H\cap K) - 1 & = & d_4(\Gamma_{H\cap K}) + d_3(\Gamma_{H\cap K})/2  - d_1(\Gamma_{H\cap K})/2\\
         & \leqslant & d_4(\Gamma_{H\cap K}) + d_3(\Gamma_{H\cap K})/2  \\
         & \leqslant & d_4(\Gamma_H)d_4(\Gamma_K) + \frac{1}{2} \left[ d_4(\Gamma_H) d_3(\Gamma_K) + d_3(\Gamma_H) d_4(\Gamma_K) \phantom{\sum_{x \in \{a,b\}^\pm}} \right. \\[-0.2in]
         &           & \phantom{d_4(\Gamma_H)d_4(\Gamma_K) + \frac{1}{2} [d_4} \left. +\!\!\!\!\!\sum_{x \in \{a,b\}^\pm} C_x(\Gamma_H) C_x(\Gamma_K) \right]\\
         & \leqslant & d_4(\Gamma_H)d_4(\Gamma_K) + d_4(\Gamma_H) d_3(\Gamma_K)/2 +
         d_3(\Gamma_H) d_4(\Gamma_K)/2 \\
         & & \phantom{d_4(\Gamma_H)d_4(\Gamma_K) + \frac{1}{2}[d_4}  + d_3(\Gamma_H) d_3(\Gamma_K)/4\\
         & \leqslant & [\rank(H) - 1] [\rank(K) - 1].
\end{eqnarray*}
and thus the conjecture holds.
\end{remark}

\section{Results}
\label{results}

\begin{remark}
\label{folding-remark} Given an endomorphism $\phi: F_2 \rightarrow
F_2$, the folding $\Gamma_{H^\phi}$ can be obtained from $\Gamma_H$ as
follows. First construct a labelled directed graph $\phi(\Gamma_H)$ by
replacing each edge with  label $x$ in $\Gamma_H$ by a sequence of
edges labelled by the successive letters of $x^\phi$ (for $x=a,b$).
Then, apply the previously described folding procedure to transform the
graph $\phi(\Gamma_{H}) \rightsquigarrow \Gamma_{H^\phi}$.  One may
verify that this yields a folding which is isomorphic to the one
obtained by constructing $\Gamma_{H^\phi}$ directly from the set $\{
w_1^\phi, \cdots, w_n^\phi \}$.
\end{remark}

For example, if $\phi$ is a length-preserving automorphism, then
$\Gamma_{H^\phi}$ can be obtained from $\Gamma_H$ by replacing every
label $x$ by $x^\phi$ and changing the orientation of the edges if
necessary.

\begin{lemma}
\label{folding-path} Let $\Gamma_H$ be the folding of a subgroup $H\fg
F_2$, and $\phi: F_2 \rightarrow F_2$ an endomorphism.  If two edges
$e,f$ from $\phi(\Gamma_{H})$ get folded during the folding process
$\phi(\Gamma_{H}) \rightsquigarrow \Gamma_{H^\phi}$, then there must
exist a path $p$ in $\phi(\Gamma_{H})$ beginning at $e$ and ending at
$f$ with the property that every edge in $p$ was folded during the
folding process.
\end{lemma}
\begin{proof}
The lemma is proved by induction on the number $n$ of edge-foldings
which take place during the folding process---note that this number
does not depend on the folding process since it is equal to
$|E_{\phi(\Gamma_{H})}| - |E_{\Gamma_{H^\phi}}|$ and the resultant
folded graph $\Gamma_{H^\phi}$ is unique).  For $n=1$, the path $p$
consists of just edges $e,f$.  Now suppose the first edge-folding
occurs when edges $d_1$ and $d_2$ are merged into an edge $d'$, and
denote the folding obtained after this identification as $\Gamma'$. By
induction, there exists a path $p'$ in $\Gamma'$ connecting $e$ and
$f$.  There are two cases to consider: either $d'$ appears in $p'$, or
it does not.  In the first case, let $p_1$ (resp. $p_2$) be the path
obtained by replacing $d'$ with $d_1$ (resp $d_2$) in $p'$.  It is
clear that either $p_1$ or $p_2$ must fulfill the requirements of the
lemma.  In the second case, we simply take $p=p'$.
\end{proof}

\begin{lemma}
\label{length-preserving-lemma}
Let $H,K\fg F_2$ be subgroups which satisfy $\delta(H,K) >
\frac{1}{2}$ (and hence are a potential counterexample to the Hanna
Neumann conjecture) and take $\ast$, $\circ$ to be two
length-preserving automorphisms of $F_2$ whose values differ on
$\mu(H,K)$, i.e. $\mu(H,K)^\ast \neq \mu(H,K)^\circ$.  Then the groups
$H^\circ, K^\ast$ must satisfy the Hanna Neumann conjecture.
\end{lemma}
\begin{proof}
Set $x_0 = \mu(H,K)$.  Since $\delta(H,K) > \frac{1}{2}$, it follows
that
\begin{eqnarray*}
  C_{x_0}(\Gamma_H) & > & \frac{1}{2} d_3(\Gamma_H) \\
  C_{x_0}(\Gamma_K) & > & \frac{1}{2} d_3(\Gamma_K).
\end{eqnarray*}

In light of Remark \ref{folding-remark}, $\Gamma_{K^\ast}$ is the
same graph as $\Gamma_K$, except that all $a$ edges have been
relabelled as $a^\ast$, and $b$ edges have been relabelled as
$b^\ast$, and an analogous statement is true about the
relationship between $\Gamma_{H^\circ}$ and $\Gamma_H$. So
\begin{eqnarray*}
  C_x(\Gamma_K) & = & C_{x^\ast}(\Gamma_{K^\ast}) \\
  C_x(\Gamma_H) & = & C_{x^\circ}(\Gamma_{H^\circ})
\end{eqnarray*}
for all $x \in \{a^{\pm},b^{\pm} \}$.  It follows that
\begin{eqnarray}
\label{x0-Kstar}
C_{x_0^\ast}(\Gamma_{K^\ast}) = C_{x_0}(\Gamma_K) > \frac{1}{2} d_3(\Gamma_K) = \frac{1}{2} d_3(\Gamma_{K^\ast}),\\
C_{x_0^\circ}(\Gamma_{H^\circ}) = C_{x_0}(\Gamma_H) > \frac{1}{2} d_3(\Gamma_H) = \frac{1}{2} d_3(\Gamma_{H^\circ}).
\end{eqnarray}

Since $x_0^\circ \neq x_0^\ast$, it follows that
$\delta_x(H^\circ, K^\ast) < \frac{1}{2}$ for every $x \in
\{a,b\}^\pm$. Thus, $\delta(H^\circ, K^\ast) < \frac{1}{2}$, and
hence by Remark \ref{walter-neumann} the groups $H^\circ,K^\ast$
cannot be a counterexample to the conjecture.
\end{proof}

\vspace{0.2in}

The proofs of {\bf Theorems \ref{first-theorem} and
\ref{first-theorem-part2}} (see page \pageref{first-theorem}) now
follow from Lemma \ref{length-preserving-lemma}.

\begin{proof}
{\em (Theorem \ref{first-theorem})}  Suppose $H,K$ do not satisfy the
conjecture.  By remark~\ref{walter-neumann} we have $\delta(H,K) >
\frac{1}{2}$.  By definition of $\tau$ we have $\delta(H^\tau,K) =
\delta(H,K)$.  We apply Lemma \ref{length-preserving-lemma} to
$H^\tau, K$, taking $\circ$ to be the fixed-point-free
length-preserving automorphism $\sigma$, and $\ast$ to be the identity
automorphism.  The theorem follows.
\end{proof}

\begin{proof}
{\em (Theorem \ref{first-theorem-part2})}  Suppose $X = \{a_1,
\ldots, a_n\}$.  We consider the embedding $\psi: F(X) \rightarrow
F_2=F(a_1,a_2)$ defined by $\psi: a_i \mapsto a_1^i a_2 a_1^i$.  If
$H, K \fg F(X)$ are a counterexample to the conjecture, then so are
$H^\psi,K^\psi \fg F_2$.  Let ${}^{-}$ be the automorphism of $F(X)$
given by $a_i \mapsto a_i^{-1}$ (for each $a_i \in X$).  Restricting
${}^{-}$ to $F_2$ and applying the previous lemma, we see that either
$H^\psi,K^\psi$ or $(H^\psi)^{-},K^\psi$ must satisfy the conjecture.
But $(H^{-})^\psi = (H^\psi)^{-}$; here we think of $H^\psi$ as a
subgroup of $F(X)$ under the canonical inclusion of $F_2$ into $F(X)$.
It follows that either $H^\psi,K^\psi$ or $(H^{-})^\psi, K^\psi$ must
satisfy the conjecture.  Since $\psi$ is a monomorphism, this implies
that either $H,K$ or $H^{-}, K$ must satisfy the conjecture.
\end{proof}

\vspace{0.2in}

Now towards the proof of Theorem \ref{second-theorem}, we introduce
the following definition:
\begin{definition}
We say $\phi: F_2(a,b) \rightarrow F_2(a,b)$ is a {\em N}-endomorphism
if it has the property that $U_\phi = \{a^\phi, b^\phi\}$ is N-reduced
\cite[pp.6]{LyndonSchupp}, which is to say that every triple $v_1$,
$v_2$, $v_3$ in $U_\phi^\pm$ satisfies
\begin{enumerate}
\item[(N0)] $v_1 \neq 1$,
\item[(N1)] $v_1 v_2 \neq 1$ implies $|v_1 v_2| \geqslant |v_1|,|v_2|$,
\item[(N2)] $v_1 v_2, v_2 v_3 \neq 1$ implies $|v_1 v_2 v_3| > |v_1| - |v_2| + |v_3|$.
\end{enumerate}
\end{definition}

\begin{remark}
\label{cannot-kill} It is well-known \cite[pp.7]{LyndonSchupp} that if
subset $U$ of a free group $F$ satisfies N0-N2, then one may associate
with each $u \in U$ words $a(u), m(u) \in F$ with $m(u) \neq 1$ such
that $u = a(u) m(u) a(u^{-1})^{-1}$ in $F$  and having the property
that for any $w = u_1 \cdots u_t$, $t \geqslant 0$, $u_i \in U^\pm$
where $u_i u_{i+1} \neq 1$, the subwords $m(u_1), \ldots, m(u_t)$
remain uncancelled in the reduced form of $w$.
\end{remark}

\begin{lemma}
\label{rank-preserving}
Every {\em N}-endomorphism of $F_2$ is a monomorphism.
\end{lemma}
\begin{proof}
Take $w \in F_2$, with $w \neq 1$.  By (N2), $| (w^\phi)^3 | > 0$,
hence $(w^\phi)^3 \neq 1$.  It follows that $w^\phi \neq 1$.
\end{proof}

\begin{lemma}
\label{small-cancellation-folding} Given $H \fg F_2$ and an {\em
N}-endomorphism $\phi$ of $F_2$, then for every edge $e$ in
$\Gamma_H$, at least one edge from the image of $e$ under $\phi$ does
not get folded during the folding process $\phi(\Gamma_{H})
\rightsquigarrow \Gamma_{H^\phi}$.
\end{lemma}
\begin{proof}
Let $e=(u,v)$ be any edge of $\Gamma_H$; suppose $e$ is labelled by
$x\in \{a, b\}^\pm$.  Consider the path $\phi(e)$ in
$\phi(\Gamma_{H})$; this path consists of a sequence of edges labelled
by successive letters of $x^\phi$.  Since $\phi$ is an {\em
N}-endomorphism, $\{a^\phi, b^\phi\}$ is N-reduced and Remark
\ref{cannot-kill} applies.  Accordingly, let $\bar{e}$ be the edge in
$\phi(e)$ which corresponds to the first letter of $m(x^\phi)$ inside
$x^\phi$.  We claim that $\bar{e}$ does not get folded during the
folding process $\phi(\Gamma_{H}) \rightsquigarrow \Gamma_{H^\phi}$.

Suppose towards contradiction, that $\bar{e}$ gets folded with some
edge $f$ during the folding process $\phi(\Gamma_{H}) \rightsquigarrow
\Gamma_{H^\phi}$.  Then by Lemma \ref{folding-path}, there must exist
a non-backtracking path $p$ in $\phi(\Gamma_{H})$ beginning at
$\bar{e}$ and ending at $f$, with the property that every edge in $p$
was folded during the folding process.  Since $p$ is a
non-backtracking path in $\phi(\Gamma_{H})$, it is a subpath of
$\phi(q)$ for some non-backtracking path $q$ in $\Gamma_{H}$.  It
follows that the labels along $\phi(q)$ are a word $u_1 \cdots u_t$,
$t \geqslant 0$, $u_i \in \{a^\phi, b^\phi\}^\pm$ and $u_i u_{i+1}
\neq 1$.  Since $\bar{e}$ is labelled by the first letter of
$m(x^\phi)$, by Remark \ref{cannot-kill} the edge $\bar{e}$ was {\em
not} folded during the folding process; this is a contradiction.
\end{proof}

We introduce the following notations: let $\Gamma_H=(V_H,E_H)$ be
the folding of $H$. Take any vertex $v\in V_H$, and let $E_v$ be
the edges incident to $v$. Define $\Gamma_v$ to be the tree
subgraph of $\Gamma_H$ induced by edges $E_v$. Then
$\phi(\Gamma_{v})$ is also a tree.

By Lemma \ref{small-cancellation-folding}, we may associate to
each edge $e\in E_v$, an edge $m(e) \in E_{\phi(\Gamma_{H})}$ which
does not get folded during the folding process $\phi(\Gamma_{H})
\rightsquigarrow \Gamma_{H^\phi}$.  We define $\text{tr}\phi(\Gamma_v)$
to be the graph obtained by truncating the branches of
$\phi(\Gamma_{v})$ so that they terminate with edges $m(e)$, $e
\in E_v$.  Clearly, for all $v \in V_H$, $\text{tr}\phi(\Gamma_v)$ is a
subgraph of $\phi(\Gamma_{H})$.

The next lemma shows that {\em N}-endomorphisms do not cause
large-scale disturbances in the neighborhood of branch vertices.

\begin{lemma}
\label{isolated-balls}
Given $H \fg F_2$ and an {\em N}-endomorphism $\phi$ of $F_2$,
then during the folding process $\phi(\Gamma_{H}) \rightsquigarrow
\Gamma_{H^\phi}$, no edge from $\text{tr}\phi(\Gamma_v)$ gets folded with
an edge from outside $\text{tr}\phi(\Gamma_v)$.
\end{lemma}
\begin{proof}
Suppose, towards contradiction, that an edge $e$ inside
$\text{tr}\phi(\Gamma_v)$ and an edge $f$ outside
$\text{tr}\phi(\Gamma_v)$ get folded during the folding process
$\phi(\Gamma_{H}) \rightsquigarrow \Gamma_{H^\phi}$.  Then by Lemma
\ref{folding-path}, there exists a path $p$ beginning at $e$ and
ending at $f$ with the property that every edge in $p$ was folded
during the folding process.  Then $p$ must pass through some edge
$m(e)$, $e \in E_v$.  This contradicts the properties of $m(e)$ as
determined in Lemma \ref{small-cancellation-folding}.  It follows that
no edge inside $\text{tr}\phi(\Gamma_v)$ gets folded with an edge
outside $\text{tr}\phi(\Gamma_v)$.
\end{proof}

Informally stated, the previous lemma implies that for an {\em
N}-endomorphism $\phi$ and subgroup $H \fg F_2$, the 5-tuple of values
\begin{eqnarray*}
  C_a(\Gamma_{H^\phi}),\; C_b(\Gamma_{H^\phi}),\; C_{a^{-1}}(\Gamma_{H^\phi}),\; C_{b^{-1}}(\Gamma_{H^\phi}),\; d_4(\Gamma_{H^\phi})
\end{eqnarray*}
is completely determinable from the 5-tuple of values
\begin{eqnarray*}
  C_a(\Gamma_H),\; C_b(\Gamma_H),\; C_{a^{-1}}(\Gamma_H),\; C_{b^{-1}}(\Gamma_H),\; d_4(\Gamma_H),
\end{eqnarray*}
without knowledge of any further structure (e.g. the generating set)
of $H$.

\vspace{0.2in}

\begin{lemma}
\label{phi0-is-magic}
Let ${\phi_0}: F_2 \rightarrow F_2$ be the endomorphism defined by
${\phi_0} a = a^2$ and ${\phi_0} b = [a,b]$.  Then for any finitely
generated subgroup $H \fg F_2$, $C_{b^{-1}}(\Gamma_{H^{\phi_0}}) =
d_3(\Gamma_{H^{\phi_0}})$, and $\text{rank}(H^{\phi_0}) =
\text{rank}(H)$.
\end{lemma}
\begin{proof}
It is straightforward to check that $\{a^2, [a,b]\}$ is N-reduced, and
hence ${\phi_0}$ is an {\em N}-endomorphism.  By Lemma
\ref{isolated-balls} it suffices to consider the effect of $\phi_0$ on
the various types of branch vertices in $\Gamma_H$.
\begin{figure}[htbp]
\centering{\mbox{\psfig{figure=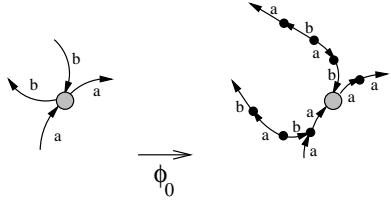}}}
\caption{The effect of $\phi_0$ on a vertex of degree $4$.}
\label{phiaction}
\end{figure}
Figure \ref{phiaction} depicts how $\phi_0$ transforms the
neighborhood of degree $4$ vertex $v$ to produce two vertices of type
$C_{b^{-1}}$.  The effect of $\phi_0$ on each of the four types of
degree $3$ vertices may be determined simply by restricting our
consideration to appropriate subgraphs of the depicted neighborhood.
It follows that $\phi_0$ transforms any degree $3$ vertex in
$\Gamma_H$ into a vertex of type $C_{b^{-1}}$ in
$\Gamma_{H^{\phi_0}}$.  Thus,
\begin{eqnarray*}
   C_{b^{-1}}(\Gamma_{H^{\phi_0}}) & = & d_3(\Gamma_H) + 2 d_4(\Gamma_H).
\end{eqnarray*}
Since all branch vertices of $\Gamma_H$ are seen to produce vertices
of type $C_{b^{-1}}$, we get that $C_{b^{-1}}(\Gamma_{H^{\phi_0}}) =
d_3(\Gamma_{H^{\phi_0}})$.  Finally, by Lemma \ref{rank-preserving},
$\phi_0$ is a monomorphism, so $\text{rank}(H^{\phi_0}) =
\text{rank}(H)$.
\end{proof}

The proof of {\bf Theorem \ref{second-theorem}} (see page
\pageref{second-theorem}) now follows from Lemma \ref{phi0-is-magic}.

\begin{proof}
{\em (Theorem \ref{second-theorem})}  Fix the embedding $\psi : F(X)
\rightarrow F_2=F(a,b)$ defined by $\psi: a_i \mapsto a^i b a^i$.  Put
$H' = (H^\psi)^{\phi_0}$ and $K' = (K^\psi)^{\phi_0}$, where
${\phi_0}$ is the endomorphism of $F(a,b)$ defined by $a^{\phi_0} =
a^2$ and $b^{\phi_0} = [a,b]$.  Then by Lemma \ref{phi0-is-magic},
\begin{eqnarray*}
C_{b^{-1}}(\Gamma_{H'}) = \frac{1}{2} d_3(\Gamma_{H'}) \\
C_{b^{-1}}(\Gamma_{K'}) = \frac{1}{2} d_3(\Gamma_{K'}),
\end{eqnarray*}
and hence $\delta(H',K')=1$.
\end{proof}

Corollary \ref{third-theorem} follows immediately since Since
$\psi\phi_0$ is a monomorphism, and hence $\text{rank}(H) =
\text{rank}(H')$, $\text{rank}(K) = \text{rank}(K')$.
\vspace{0.2in}

\section{Acknowledgements}

The first two authors are grateful to the Mathematics department at
the City University of New York Graduate Center for funding this work
as part of the their ongoing doctoral research.  The second author
thanks both the Center for Computational Sciences at the Naval
Research Laboratory in Washington DC, and the Advanced Engineering and
Sciences division of ITT Industries for their support of these
endeavors.

\bibliographystyle{plain}
\bibliography{hnc-JKM}
\end{document}